\documentclass[12pt,a4paper]{article}

\usepackage{amsmath,amscd,amssymb,amsthm,amstext,mathrsfs}
\input{xypic.tex}

\theoremstyle{plain}
   \newtheorem{theorem}{Theorem}[section]
   
   \newtheorem{lemma}[theorem]{Lemma}
   
   \theoremstyle{definition}
   
   \newtheorem{definition}[theorem]{Definition}
   
   \theoremstyle{remark}
   
   \newtheorem{remark}[theorem]{Remark}

\newcommand{\FF}{{\mathbb F}}

\newcommand{\simp}{_\bullet}

\newcommand{\gshuf}{{V}}

\newcommand{\id}{{\mathrm {id}}}
\newcommand{\Id}{{\mathrm {Id}}}
\newcommand{\susp}{{\cal S}}
\newcommand{\opsusp}{{\cal S}^\prime }
\newcommand{\tw}{{\cal T}}
\newcommand{\EZtype}{EM type transformation}

\title{An alternative approach to homotopy operations}
\author{Marcel B\" okstedt \& Iver Ottosen}
\date{May 17, 2004}

\begin{document}
\maketitle

\begin{abstract}
We give a particular choice of the higher
Eilenberg-MacLane maps by a recursive formula.
This choice leads to a simple description of the
homotopy operations for simplicial ${\bf Z}/2$-algebras.
\end{abstract}

\section{Introduction}
This paper is about the ring of homotopy groups of a simplicial ring.
This ring of homotopy classes has a huge amount of additional structure. 
The theory is best worked out for algebras over $\FF_2$, and
we will restrict ourselves to this case.
\S\S 2-3 in \cite{G} contain a good survey with references 
to the original articles. We just recall the points which
are most important to us.

The main observation is that the square of every element in 
positive degree is zero. Analyzing this fact gives rise to a divided 
power structure on the ring of homotopy groups.
There is a refinement of this, which constructs 
a sequence of homotopy operations $\delta_i$. These
are defined by Dwyer in \cite{D} and also by Bousfield.

The homotopy of a simplicial algebra $R\simp$ is isomorphic to
the homology of the associated chain complex $C_*(R)$.
We can represent an element in $\pi_n (R)$ by a cycle in the
the chain complex, that is by a class
$z\in R_n=C_n(R)$, such that $\sum_{0\leq i \leq n}d_i(z)=0$.

The problem which this article wants to solve is the following:
Suppose that $\delta_i$ is defined on $\pi_n(R\simp)$.
Can we give an explicit formula for an element
in $R_{n+i}$ that represents the element
$\delta_i(z)$?

The reason that we care is that we are interested
in explicit applications, like in \cite{BO} where we compute 
these operations.

We will focus on Dwyer's approach. In order to define the 
operations, he considers a sequence of natural transformations, 
defined for pairs of simplicial $\FF_2$-vector spaces 
$V\simp$ and $W\simp$ as follows:
\[
D^k\colon [C(V)\otimes C(W)]_{m}\to [C(V\otimes W)]_{m-k}.
\]  
The transformations satisfy recursive conditions, which we write 
down in the beginning of section \ref{sec:HEM}. 

These maps are considered to be higher analogues of the
Eilenberg-Mac- Lane map. Dwyer proves an existence and
uniqueness result for these maps, but he does not
give an explicit formula. He then uses them to define 
homotopy operations. We follow the reformulation of \cite{G} 
at this point.

Let $z\in R_n$ and define $\Theta_i(z)$ (for $1\leq i\leq n$)
as the element 
\[
\Theta_i(z)=\mu D^{n-i}(z\otimes z)+\mu D^{n-i-1}(z\otimes \partial z).
\]
Here $\mu$ is multiplication in the simplicial algebra.
If $z$ is a cycle, and $2\leq i\leq n$, one sees that $\Theta_i(z)$  
is also a cycle, and that the formula defines an operation
$\delta_i\colon \pi_n(R\simp)\to \pi_{n+i}(R\simp)$.

We are going to give simpler formulas for 
$D^i$ which satisfy the defining recursive relations.
The derivation of these is where the hard work of 
this article is done. When we have obtained the formulas,
we can plug them into the definition of $\Theta_i$, and obtain 
an explicit formula for $\delta_i(z)$. 
We now explain the formula we obtain in this fashion.

Write $N_*R$ for the normalized chain complex with
$N_qR= \cap_{1\leq i \leq q} \ker (d_i)$ and 
differential $d_0$. Let $Z_qR\subseteq N_qR$ denote
the cycles i.e. the elements $z\in R_q$ with
$d_iz=0$ for $0\leq i \leq q$.

To make the final formulas appear simpler, we assume that 
$z$ is a normalized chain. This is no restriction, since
the associated chain complex considered above contains the 
normalized chain complex as a quasi-isomorphic subcomplex.

\begin{definition}
For integers $q$, $i$ with $1\leq i \leq q$ we define $U(q,i)$ to 
be the set of pairs $(\mu , \nu )$ of ordered sequences
$\mu_1 < \dots < \mu_i$, $\nu_1  < \dots < \nu_i$
with disjoint union
$$\{ \mu_1 , \dots ,\mu_i \} \sqcup \{ \nu_1 ,\dots , \nu_i \} =
\{ q-i, q-i+1, \dots ,q+i-1 \}.$$
Let $\gshuf (q,i)\subseteq U(q,i)$ be the subset with $\mu_1=q-i$.
\end{definition}

\begin{definition}
\label{deltadef}
For a cycle $z\in Z_qR$ we define $\delta_i(z) \in R_{q+i}$ by
$$\delta_i (z) = \sum_{ (\mu ,\nu ) \in \gshuf (q,i)}
s_{\nu_i} \dots s_{\nu_1} (z) s_{\mu_i} \dots s_{\mu_1} (z).$$
\end{definition}

Note the close relationship to the Eilenberg-MacLane map $D$: 
\[
\mu D(z\otimes z)= \sum_{ (\mu ,\nu ) \in U(q,q)}
s_{\nu_i} \dots s_{\nu_1} (z) s_{\mu_i} \dots s_{\mu_1} (z).
\]

\begin{theorem}
When $2\leq i \leq q$ the formula in Definition 
\ref{deltadef} defines a map $\delta_i : Z_qR \to Z_{q+i}R$. 
The induced map on homotopy 
$\delta_i : \pi_q R \to \pi_{q+i} R$ is the Bousfield-Dwyer
homotopy operation.
\end{theorem}

\begin{proof}
Let $z\in N_qR$ be a cycle. It is shown in Lemma \ref{ztoz} 
that $\delta_i(z)$ is a cycle in $N_{q+i}R$. 

We use Definition \ref{def:D} as our choice of higher 
Eilenberg-MacLane maps. The suspension operator $\susp$
in the definition increases the simplicial degree by one.
It preserves composition of simplicial maps and 
$\susp (d_j)=d_{j+1}$, $\susp (s_j)=s_{j+1}$, $\susp (\id ) = \id$. 

We have that $\Theta_i (z) = \mu D^{q-i}(z\otimes z)$. The 
statement we must show is that 
$\Theta_i (z) = \mu \susp^{q-i} (D^0) (z\otimes z)$. So it
suffices to see that $D^k(z\otimes z) = S^k(D^0)(z\otimes z)$,
where $q-i=k$. 

For $0\leq j \leq k-1$ we have 
\begin{align*}
\susp^j (D^{k-j}) &= \susp^{j+1}(D^{k-j-1}) +
\begin{cases}
\susp^j (D^{k-j-1}(d_0\otimes \id )) &, k-j \text{ even} \\
\susp^j (D^{k-j-1}(\id \otimes d_0)) &, k-j \text{ odd}
\end{cases}
\\
&= \susp^{j+1}(D^{k-j-1}) +
\begin{cases}
\susp^j (D^{k-j-1})(d_j\otimes \id ) &, k-j \text{ even} \\
\susp^j (D^{k-j-1})(\id \otimes d_j) &, k-j \text{ odd.}
\end{cases}
\end{align*}
Thus $\susp^j (D^{k-j})(z\otimes z)=
\susp^{j+1} (D^{k-j-1})(z\otimes z)$. We iterate this 
result and find that
$D^k(z\otimes z)=\susp^0 (D^k) (z\otimes z) = 
\susp^k (D^0) (z\otimes z)$.
\end{proof}

\section{Higher Eilenberg-MacLane maps.}
\label{sec:HEM}
Let $V_\bullet$, $W_\bullet $ be simplicial $\FF_2$-vector spaces.
The Eilenberg-MacLane map $D$ can be considered as a set of
linear maps
\begin{equation}
  \label{eq:shuffle}
D_{i,j}\colon V_i \otimes W_j \to V_{i+j}\otimes W_{i+j}.  
\end{equation}
There is an explicit formula for these maps
\[
D_{i,j}=\sum s_{\nu_j} \dots s_{\nu_1} \otimes s_{\mu_i}\dots s_{\mu_1},
\]
where the sum is indexed by $(i,j)$-shuffles. An $(i,j)$-shuffle
consists of two increasing sequences $\mu =(\mu_1 ,\mu_2 ,\dots \mu_i )$ 
and $\nu =(\nu_1 ,\nu_2 ,\dots ,\nu_j )$ such that each integer from
$\{ 0,1,\dots ,i+j-1\}$ occurs exactly once as either a $\mu_r$
or a $\nu_s$. 
There is an excellent discussion of this map in \cite{ML}, chap. VIII \S 8. 

We will also consider the map $\phi_k$ with 
$(\phi_k)_{i,j} \colon V_i\otimes W_j \to V_{i+j-k}\otimes W_{i+j-k}$,
defined to be the zero map, unless $i=j=k$. In this exceptional case,
the map is the identity. 

Let $T_{i,j}\colon V_i\otimes W_j \to W_j \otimes V_i$ 
be the map permuting the factors. According to
\cite{D} there is a sequence of maps $D^k$, $k=0,1,2,\dots $ with
\[
D^k_{i,j}\colon V_i\otimes W_j \to V_{i+j-k}\otimes W_{i+j-k},
\] 
only defined under the conditions that
\begin{equation}
\label{eq:dimensions}
0 \leq 2k \leq i+j,
\end{equation}
which satisfy that
\begin{equation}
\label{eq:Dwyercondition0}
D^0+TD^0T=D+\phi_0,
\end{equation}
\begin{equation}
\label{eq:Dwyerconditions}
D^k+TD^kT=D^{k-1}\partial+\partial D^{k-1}+\phi_k.
\end{equation}

There are many choices for the maps $D^k$, but Dwyer proves that
the choices are equivalent up to homotopy (in a strong sense). 


Before we give our definition of the maps $D^k$, we introduce
some language, for the purpose of avoiding a nightmare of indices.
We are going to write down various formal sums of products of the
simplicial generators $d_r,s_r$. Such an expression
sometimes but not always defines a map $V_i\to V_j$
for all simplicial vector spaces $V\simp$.
For instance, consider the simplicial relation $d_rs_r=\id$.
The right hand side of this relation defines a map (the identity) 
$V_i \to V_i$ for all $i$. The left hand side does 
not define a map on $V_i$ if $r>i$.
But if $r\leq i$, both sides of the equation defines maps, and
in this case the relation says that the 
two maps $V_i \to V_i$ induced by the two sides agree. 

We are going to consider tensor products of pairs of simplicial groups
$V\simp\otimes W\simp$.
For pairs of integers $(i,j),(k,l)$, we will write
down natural transformations 
$V_i\otimes W_j\to V_k\otimes W_l$
and relations between such. We will specify them 
as sums of formal sequences of generators $d_r,s_r$. 
Every time we do this, we have to keep
track of whether the formal sequences do indeed define natural 
transformations of the groups we write down.

We start by introducing a ``suspension operator'' $\susp$.
\begin{definition}
\label{def:susp}
Let $\Delta$ denote the simplicial category.
Define the functor $\opsusp :\Delta \to \Delta$ by 
$\opsusp([n])=[n+1]$ on objects, and 
\[
(\opsusp \alpha )(i)=
\begin{cases}
\alpha (i-1)+1 &\text{ if } i\geq 1 \\
0 &\text{ if } i=0. 
\end{cases}
\]  
on morphisms. Let $\susp = (\opsusp )^{op}:\Delta^{op} \to
\Delta^{op}$ be the corresponding functor on the opposite category.
\end{definition}

We can picture the suspension operator as follows:
\vskip 12pt
\hfil
\raisebox{-0.7cm}
{\begin{minipage}{5cm}
\[
\diagram
2 \drto & 2\\
1 \rto  & 1\\
0 \urto & 0\\
\enddiagram
\]

\[
\alpha \colon[2]\to [2]
\]
\end{minipage}
}
\hfil
\begin{minipage}{5cm}
\[
\diagram
3 \drto & 3\\
2 \rto  & 2\\
1 \urto & 1\\
0 \rto  & 0
\enddiagram
\]

\[
\opsusp \alpha \colon[3]\to [3]
\]
\end{minipage}
\hfil
\vskip 12pt

Note that $d^0$ defines a natural transformation 
$d^0: \Id \to \opsusp$. So we have a natural transformation
$d_0 : \susp \to \Id$.

For a simplicial vector space 
$V \colon \Delta^{op} \to \{ \text{$\FF_2$-vector spaces} \}$, 
we define a new simplicial vector space $\susp (V)\simp $ 
as the functor  $V\circ \susp$. There is
a natural transformation of simplicial vector spaces
$P\colon \susp (V)\simp \to V\simp$, given 
degreewise as $P\colon \susp (V)_i=V_{i+1} \overset{d_0}\to V_i$.

Consider any natural transformation defined on the category of
simplicial vector spaces of the form
\[
\theta_V \colon V_i \to V_j,
\]
for instance, this could be a map 
induced by a morphism $[j]\to [i]$ in $\Delta$, or a linear combination of
such maps.
\begin{definition}
\label{def:suspension}
The suspension of $\theta$ is the natural transformation
\[
(\susp \theta )_V  \colon V_{i+1} = \susp (V)_{i}
\overset {\theta_{\susp V}}\longrightarrow 
\susp (V)_j = V_{j+1}. 
\] 
\end{definition}

We obtain a commutative diagram:
\[
\begin{CD}
V_{i+1} @= \susp (V)_i @>{P}>> V_i \\
@V{(\susp \theta )_V}VV @V{\theta_{\susp V}} VV @V{\theta_V}VV \\
V_{j+1} @= \susp (V)_j @>{P}>> V_j \\
\end{CD}
\]
The composite of both rows is $d_0$, so it follows from this diagram 
and the naturality of $\theta$ that we have a relation
\begin{equation}
\label{eq:d0andtheta}
d_0(\susp\theta)=\theta d_0\colon \susp (V)_{i}\to V_j.  
\end{equation}

We can also define a suspension on {\em formal} products of
the simplicial generators $d_i,s_i$ by $\susp(d_i)=d_{i+1}$ 
and $\susp (s_i)=s_{i+1}$. This suspension is compatible with the 
(formal) simplicial relations. In this way, we can
formally define $\susp$ on strings of simplicial generators.
by adding 1 to the index of every occurring generator $s_i$
or $d_j$.
   
\begin{lemma}
Let $V\simp$ be a simplicial $\FF_2$-vector space.
\begin{itemize}
\item Let $\alpha$ be a formal products of generators $d_r,s_r$.    
Then $\alpha$ defines a natural transformation 
$\alpha_* \colon V_i\to V_j$ if and only if
the corresponding natural transformation 
$\susp (\alpha )_* :V_{i+1}\to V_{j+1}$ is also defined.
\item Suspension commutes with passage to natural 
transformation, if the natural transformation is defined. 
That is
\[
\susp (\alpha_* )=(\susp \alpha)_*.
\]
\item If a formal sum $\sum \alpha_* \colon V_i\to V_j$  
is defined and equals the zero map for $i$,
then $\sum \susp \alpha_* \colon V_{i+1}\to V_{j+1}$ 
is defined, and also equals the zero map.
\end{itemize}
\end{lemma}

\begin{proof}
By induction on the number of factors in a product of generators, 
it is enough to check the first statement for the case of a generator 
$s_j$ or $d_j$. 
For example, $(s_j)_*$ is defined on $V_k$ if and only if
$k \geq j$. But it is also true that  
$(\susp (s_j))_*=(s_{j+1})_*$ is defined on $V_{k+1}$ if and only if
$k+1\geq j+1$. 
The case $d_j$ is treated in exactly the same way.

To prove the second statement, by induction it is enough
to consider the case of a simplicial generator.
That is, we have to check that $\susp (s_i)=s_{i+1}$ and that
$\susp (d_i)=d_{i+1}$. But this follows directly from the
definition of $\susp$.

The third statement follows from the second, since
\[
(\sum \susp \alpha)_*=\susp (\sum \alpha_* )=\susp 0=0
\]
\end{proof}

All we have said can be generalized to tensor products of
two simplicial vector spaces.

\begin{lemma}
\label{lemma:zeroes}
Let $V\simp$, $W\simp$ be simplicial $\FF_2$-vector spaces.
We can define the suspension of a natural transformation
$\theta \colon V_i\otimes W_j \to V_k \otimes W_l$ as a
natural transformation
$\susp \theta \colon V_{i+1}\otimes W_{j+1} \to V_{k+1} \otimes W_{l+1}$.
\begin{itemize}
\item Let $\alpha , \beta$ be formal products of generators $d_i,s_i$.    
If $\alpha_* \otimes \beta_*$ is defined for the index
$(i,j)$ then $\susp (\alpha_* \otimes \beta_* )$
is defined for $(i+1,j+1)$.
\item Suspension commutes with passage to natural transformation.
\item If a formal sum $\sum \alpha_* \otimes \beta_*$  
is defined and equals the zero map for $(i,j)$,
then $\susp (\sum \alpha_* \otimes \beta_* )$ equals the zero map
for $(i+1, j+1)$.
\end{itemize}
\end{lemma}
We leave the proof to the dedicated reader.

Note that the map $(\phi_k)_{i,j}$ considered above is a natural
transformation which is not defined as a formal combination of
the simplicial generators $d_i,s_i$. But the suspension on it is 
still defined, and actually
\[
(\phi_{k+1} )_{i+1,j+1}=\susp ((\phi_k )_{i,j}) \colon 
V_{i+1}\otimes W_{j+1} \to V_{i+j+1-k} \otimes W_{i+j+1-k}.
\]

However, $\phi_k$ is an example of a map with the following property.
\begin{definition}
An \EZtype{} $F$ consists of the following data.
\begin{enumerate}
\item An index function $I_F$ which to each pair of integers $(i,j)$ 
associates a pair of integers $(k,l)=I_F(i,j)$.
\item For each $(i,j)$, we have a natural transformation on the category
of simplicial vector spaces $V_i\otimes W_j \to V_k\otimes W_l$.
Here, we conventionally define $V_i\otimes W_j=0$ if either $i<0$ or $j<0$. 
\end{enumerate}
\end{definition}

For instance, $\phi_k$ defines such an \EZtype . 
The index function for $\phi_k$ is $I_{\phi_k}(i,j)= (i+j-k,i+j-k)$.

The main example of such a map  is the Eilenberg-MacLane map $D$.
It is a collection of maps 
$D_{i,j} \colon V_i\otimes W_j \to V_{i+j}\otimes W_{i+j}$. 
The index function is $I_D(i,j)=(i+j,i+j)$.

We can always compose two \EZtype s $F,G$.
We can add them if the index function of $F$ agrees
with the index function of $G$. If $F$ and $G$ does
not have the same index function, their sum
is not defined. This is the price we pay for 
keeping easy control of the indices involved.
 
We can suspend an \EZtype. We define
$I_{\susp F}(i+1,j+1)=(1,1)+I_F(i,j)$ and 
$(\susp F)_{i+1,j+1}=\susp(F_{i,j})$. This is to
be interpreted so that if either $i<0$ or
$j< 0$, then $(\susp F)_{i+1,j+1}=0$.

We can also twist an \EZtype{} by defining $\tw (F)_{i,j}=TF_{j,i}T$.
Suspension and twisting preserve sum and composition, 
and they commute.
 
Here are some examples of \EZtype s, and relations between them. 
We insist that the relations are valid as relations between
natural transformations with source  $V_i\otimes W_j$ for all 
pairs of integers $(i,j)$.
When checking the formulas below, the main worry is 
to keep track of cases like $V_0\otimes W_j$ and $V_i\otimes W_0$

The sequences $d_0\otimes \id$,
$s_0\otimes \id$ are \EZtype s, with index functions
$(i,j)\mapsto (i-1,j)$ respectively $(i,j)\mapsto (i+1,j)$.  
The thing to notice is that
$d_0\otimes \id \colon V_0\otimes W_j \to V_{-1}\otimes W_j$,
makes sense (and equals the zero map) 
since we define $V_i=0$ for $i=-1$, and
$W_j=0$ for $j=-1$.  

The simplicial relation give relations between these functors. 
\begin{equation}
\label{eq:simp0}
(d_0\otimes \id )(s_0\otimes \id )=\id \otimes \id
\end{equation}
is true as a relation between \EZtype s.

Another type of example is $\partial \otimes \id$,
with index function $(i,j)\mapsto (i-1,j)$. It is defined by
$(\partial\otimes \id )_{i,j}=\sum_{0\leq r\leq i} d_r\otimes \id$.
Similarly, we can define $\id \otimes \partial$, with
index function $(i,j)\mapsto (i,j-1)$.

When we apply suspension, we have to be
careful. Here is an example of this: {\em unless} $i>0$, $j=0$,
we have that
$\susp (\partial \otimes \id )_{i,j} + (d_0\otimes \id )_{i,j}=
(\partial \otimes \id )_{i,j}$. If we post-compose with 
an \EZtype{} which vanishes on all groups $V_i\otimes W_0$, we get a
genuine relation. For instance we have for any \EZtype{} $F$ that:  
\begin{equation}
\label{eq:simp1}
\susp F\susp (\partial \otimes \id ) + \susp F(d_0\otimes \id )=
\susp F(\partial \otimes \id ).
\end{equation}

We can also use simplicial operations simultaneously in both factors.
Using the index function $(i,j)\mapsto (i-1,j-1)$, we put 
\[
\delta = \sum_{0\leq r\leq \min(i,j)} {d_r\otimes d_r \colon
V_i\otimes W_j \to V_{i-1}\otimes W_{j-1}}.  
\]
$d_0\otimes d_0$ is another \EZtype{} with the same index function, and
\begin{equation}
\label{eq:simp4}
\susp (\delta )=\delta +d_0\otimes d_0.
\end{equation}

Similarly (using the simplicial relations $d_0s_0=\id =d_1s_0$,
$d_is_0=s_0d_{i-1}$ for $i\geq 2$) we get that
\begin{equation}
\label{eq:simp3}
(\id \otimes \partial)(\id \otimes s_0)=
(\id \otimes s_0)(\id \otimes \partial)+
(\id \otimes s_0d_0).
\end{equation}
In the same way
\begin{equation}
  \label{eq:simp2}
  (\id\otimes \partial)(\id\otimes d_0)=
(\id\otimes d_0)(\id\otimes \partial)+
(\id\otimes d_0d_0).
\end{equation}
For any \EZtype{} $F$, we have 
(because of (\ref{eq:d0andtheta})) 
\begin{equation}
  \label{eq:simp5}
  (d_0\otimes d_0)(\susp F)=F(d_0\otimes d_0).
\end{equation}

We now give our definition of the higher Eilenberg-MacLane maps.
\begin{definition}
\label{def:D}
$D^k$ is the \EZtype{} with index function
$I_{D^k}(i,j)=(i+j-k,i+j-k)$, and defined
by $D^0=\susp (D) (\id \otimes s_0)$ and
inductively for $k\geq 1$ by the formula
\begin{equation}
\label{definition}
D^k=\susp (D^{k-1})+
\begin{cases}
D^{k-1}(d_0\otimes \id ) &\text{if $k$ is even} \\
D^{k-1}(\id \otimes d_0) &\text{if $k$ is odd.} \\  
\end{cases}
\end{equation}  
\end{definition}

The sum on the right hand side of (\ref{definition}) 
is defined because the \EZtype s
$\susp (D^{k-1})$ and $D^{k-1}(\id\otimes d_0)$ have the same
index function.

This equation (\ref{definition}) is an equation of \EZtype s.
If we write it out, it means that we have 
natural transformations $D^k_{i,j}\colon V_i\otimes W_j \to V_{i+j-k}\otimes
W_{i+j-k}$
given inductively as
\begin{equation*}
D^k_{i,j}=
\begin{cases}
\susp (D^{k-1}_{i-1,j-1})+D^{k-1}_{i-1,j}(d_0\otimes \id)&\text{if $k$ is even, $i,j\geq 1$,}\\
\susp (D^{k-1}_{i-1,j-1})+D^{k-1}_{i,j-1}(\id\otimes d_0)&\text{if $k$ is odd, $i,j\geq 1$,}\\
D^{k-1}_{i-1,j}(d_0\otimes \id)&\text{if $k$ is even, $i\geq 1,j=0$,}\\
D^{k-1}_{i,j-1}(\id\otimes d_0)&\text{if $k$ is odd, $i=0,j\geq 1$,}\\
0&\text{else}.  
\end{cases}
\end{equation*}

It is obvious that our $D^0$ satisfies (\ref{eq:Dwyercondition0}).
We have to prove that our choices of $D^k$, $k\geq 1$ satisfy 
(\ref{eq:Dwyerconditions}). Our strategy for proving this is first 
to prove relations between \EZtype s, that is a relation between 
natural transformations valid for all pair of integers $(i,j)$.

\begin{definition}
Let $A^k$ be the \EZtype{} with index function 
$I_{A^k}(i,j)=(i+j-k,i+j-k)$, defined by
\begin{equation*}
A^0=D^0+\tw D^0+D,
\end{equation*}
\begin{equation*}
A^k=D^k+\tw D^k +\delta D^{k-1}+
D^{k-1}(\partial \otimes \id )+D^{k-1}(\id \otimes \partial ).
\end{equation*}
\end{definition}

We have now defined all \EZtype s that we will need.
The main work of this article is done in proving the 
following recursion relation for $A^k$:

\begin{lemma}
\label{lemma:formula}
For $k\geq 1$ we have that
\begin{equation*}
A^k=\susp (A^{k-1})+
\begin{cases}
A^{k-1}(\id \otimes d_0) & \text{ if $k$ is even} \\
A^{k-1}(d_0\otimes \id ) & \text{ if $k$ is odd.} \\  
\end{cases}
\end{equation*}
\end{lemma}

\begin{proof}
The proof is by direct computation. It is divided into three
cases: $k=1$; $k\geq 2$ and $k$ even; $k\geq 3$ and $k$ odd.
The method used in the three cases is the same.
We will write down about a dozen relations, and then
add all of them to give the desired recursion formula.

{\it The case $k=1$.}   
We want to prove that 
\[
\tag{*}
\begin{split}
& D^1+\tw D^1+D^0(\partial \otimes \id )
+D^0(\id \otimes \partial)+\delta D^0 \\
&+\susp (D^0+\tw D^0+D)+(D^0+\tw D^0+D)(d_0\otimes \id )=0. 
\end{split}
\]

We have that $D^0+(\susp D)(\id \otimes s_0)=0$.
From this follows immediately the relations
\begin{equation}
\label{eq:rel1}
\delta D^0 + \delta (\susp D)(\id \otimes s_0)=0, 
\end{equation}
\begin{equation}
\label{eq:rel2}
D^0(\id \otimes d_0)+\susp D(\id \otimes s_0d_0)=0,
\end{equation}
\begin{equation}
\label{eq:rel3a}
D^0(\partial \otimes \id )+\susp D(\partial \otimes s_0)=0,
\end{equation}
\begin{equation}
\label{eq:rel3b}
D^0(\id \otimes \partial)+
\susp D(\id \otimes s_0)(\id \otimes \partial)=0,
\end{equation}
and 
\begin{equation}
\label{eq:rel4}
D^0(d_0\otimes \id )+\susp D(d_0 \otimes s_0)=0.
\end{equation}

There is also the defining relation for $D^1$, and we can apply the
twist to it. This gives
\begin{equation}
D^1+\susp D^0+D^0(\id \otimes d_0)=0,
\end{equation}
and 
\begin{equation}
\tw D^1+ \susp \tw D^0+ \tw D^0(d_0\otimes \id )=0.
\end{equation}

A fundamental property of the  Eilenberg-MacLane map is 
that it is a chain map. In our notation, this says that
$\delta D+D(\partial \otimes \id)+D(\id \otimes \partial )=0$.

Suspending this, using \ref{lemma:zeroes}, right multiplying with 
$\id\otimes s_0$ we get
\begin{equation}
\label{eq:rel7}
\susp (\delta D) (\id \otimes s_0)+
\susp (D(\partial \otimes \id ))(\id \otimes s_0)+
\susp (D(\id \otimes \partial ))(\id\otimes s_0)=0.
\end{equation}
Now use that $\susp$ is compatible with composition so that
$\susp (D(\partial \otimes \id ))=
(\susp D)\susp (\partial \otimes \id )$. 
The identitity (\ref{eq:simp1}) provides
\begin{equation}
\label{eq:rel8a}
\susp (D(\partial \otimes \id ))(\id \otimes s_0)+
\susp D(\partial \otimes s_0)+
\susp D(d_0\otimes s_0)=0.
\end{equation}
Similarly, using twisted versions of
(\ref{eq:simp1}) and (\ref{eq:simp0}) gives
\begin{equation}
\label{eq:rel8b}
\susp (D(\id \otimes \partial ))(\id \otimes s_0)+
\susp D (\id \otimes \partial )(\id \otimes s_0)+
\susp D=0.
\end{equation}
The relation (\ref{eq:simp3}) gives
\begin{equation}
\label{eq:rel9}
\susp D(\id \otimes \partial )(\id \otimes s_0)+
\susp D(\id \otimes s_0)(\id \otimes \partial)+
\susp D(\id \otimes s_0d_0)=0. 
\end{equation}

The relation (\ref{eq:simp4}) also gives a relation. We apply
(\ref{eq:simp5}) for $F=D$ to it, and obtain:
\begin{equation}
\label{eq:rel10}
\delta (\susp D)(\id \otimes s_0)+
\susp (\delta D)(\id \otimes s_0)+
D(d_0\otimes \id )=0.
\end{equation}
Adding the numbered relations (\ref{eq:rel1})-(\ref{eq:rel10})
gives (*), and finishes the proof of case $k=1$.

{\it The case $k\geq 2$ and $k$ even.} 
The relation we want to  prove is the following:
\begin{equation}
\tag{**}
\begin{split}
& D^k+\tw D^k+D^{k-1}(\partial \otimes \id )+
D^{k-1}(\id \otimes \partial )+\delta D^{k-1} \\
& +\susp (D^{k-1}+\tw D^{k-1}+D^{k-2}(\partial \otimes \id )
+D^{k-2}(\id \otimes \partial )+\delta D^{k-2}) \\
& +(D^{k-1}+\tw D^{k-1}+D^{k-2}(\partial \otimes \id )
+D^{k-2}(\id \otimes \partial )+\delta D^{k-2})(\id \otimes d_0)=0.
\end{split}
\end{equation}
The definition says that 
\begin{equation}
\label{eq:firsteven}
D^k+\susp D^{k-1}+D^{k-1}(d_0\otimes \id )=0.  
\end{equation}
Applying the twist to this, we get
\begin{equation}
\tw D^k+\susp \tw D^{k-1}+\tw D^{k-1}(\id \otimes d_0)=0.  
\end{equation}
Since $k-1$ is odd, the definition says that
$D^{k-1}=\susp D^{k-2}+D^{k-2}(\id \otimes d_0)$. This gives us
a sequence of relations:
\begin{equation}
D^{k-1}(d_0\otimes \id ) +\susp D^{k-2}(d_0\otimes \id )+
D^{k-2}(d_0\otimes d_0)=0,
\end{equation}
\begin{equation}
D^{k-1}(\id \otimes d_0)+\susp D^{k-2}(\id \otimes d_0)+
D^{k-2}(\id \otimes d_0d_0)=0,
\end{equation}
\begin{equation}
\delta D^{k-1}+\delta D^{k-2}(\id \otimes d_0)+\delta \susp D^{k-2}=0,
\end{equation}
\begin{equation}
D^{k-1}(\partial \otimes \id )+
D^{k-2}(\partial \otimes \id )(\id \otimes d_0)
+\susp D^{k-2}(\partial \otimes \id )=0,
\end{equation}
and
\begin{equation}
D^{k-1}(\id \otimes \partial )+\susp D^{k-2}(\id \otimes \partial )
+D^{k-2}(\id \otimes d_0\partial )=0.
\end{equation}
(\ref{eq:simp2}) gives us
\begin{equation}
D^{k-2}(\id \otimes \partial )(\id \otimes d_0)+
D^{k-2}(\id \otimes d_0\partial )+
D^{k-2}(\id \otimes d_0d_0)=0.
\end{equation}
Using (\ref{eq:simp5}) for $F=D^{k-2}$ we get
\begin{equation}
\label{eq:firstsymmetric}
(d_0\otimes d_0)\susp D^{k-2}+D^{k-2}(d_0\otimes d_0)=0.
\end{equation}
Finally, the following equations follow from (\ref{eq:simp1}) and
its twisted version since
in this case $k\geq 3$ and from (\ref{eq:simp4}),
using the compatibility of suspension with products.
\begin{equation} 
\susp (D^{k-2}(\partial \otimes \id ))+
\susp D^{k-2}(\partial \otimes \id )+
\susp D^{k-2}(d_0\otimes \id ),
\end{equation}
\begin{equation} 
\susp (D^{k-2}(\id \otimes \partial ))+
\susp D^{k-2}(\id \otimes \partial )+
\susp D^{k-2}(\id \otimes d_0),
\end{equation}
and
\begin{equation}
\label{eq:lasteven} 
\susp (\delta D^{k-2})+
\delta \susp D^{k-2}+
(d_0\otimes d_0)\susp D^{k-2}.
\end{equation}
If we add all numbered equations 
(\ref{eq:firsteven})-(\ref{eq:lasteven}), we get formula (**). This finishes
the case $k\geq 2$, $k$ even.

{\it Case $k\geq 3$ and $k$ odd.} In this case, we want to prove that
\begin{equation}
\tag{***}
\begin{split}
& D^k+\tw D^k+D^{k-1}(\partial \otimes \id )+
D^{k-1}(\id \otimes \partial )+\delta D^{k-1} \\
& +\susp (D^{k-1}+\tw D^{k-1}+D^{k-2}(\partial \otimes \id )
+D^{k-2}(\id \otimes \partial )+\delta D^{k-2}) \\
& +(D^{k-1}+\tw D^{k-1}+D^{k-2}(\partial \otimes \id )
+D^{k-2}(\id \otimes \partial )+\delta D^{k-2})(d_0\otimes \id )=0.
\end{split}
\end{equation}

The definition and its twist give
\begin{equation}
\label{eq:firstodd}
D^k+\susp (D^{k-1})+D^{k-1}(\id \otimes d_0)=0,
\end{equation}
\begin{equation}
\tw D^k+\susp \tw (D^{k-1})+\tw D^{k-1}(d_0\otimes \id )=0.
\end{equation}
Since $D^{k-1}=\susp(D^{k-2})+D^{k-2}(d_0\otimes \id )$, 
we have relations:
\begin{equation}
D^{k-1}(\id \otimes d_0) +\susp D^{k-2}(\id \otimes d_0)+
D^{k-2}(d_0\otimes d_0)=0,
\end{equation}
\begin{equation}
D^{k-1}(d_0\otimes \id )+\susp D^{k-2}(d_0\otimes \id )+
D^{k-2}(d_0d_0\otimes \id )=0,
\end{equation}
\begin{equation}
\delta D^{k-1}+\delta D^{k-2}(d_0\otimes \id )+\delta \susp D^{k-2}=0,
\end{equation}
\begin{equation}
D^{k-1}(\partial \otimes \id )+
D^{k-2}(d_0 \partial \otimes \id )
+\susp D^{k-2}(\partial \otimes \id )=0,
\end{equation}
and
\begin{equation}
D^{k-1}(\id \otimes \partial )+
D^{k-2}(\id \otimes \partial )(d_0\otimes \id)+
\susp D^{k-2}(\id \otimes \partial )=0.
\end{equation}
The twisted version of (\ref{eq:simp2}) gives us
\begin{equation}
\label{eq:lastodd}
D^{k-2}(\partial \otimes \id )(d_0\otimes \id )+
D^{k-2}(\partial d_0\otimes \id )+
D^{k-2}(d_0d_0 \otimes \id )=0.
\end{equation}
Adding the numbered equations (\ref{eq:firstodd})-(\ref{eq:lastodd})
and the numbered equations 
(\ref{eq:firstsymmetric})-(\ref{eq:lasteven}), we get (***)
\end{proof}

\begin{theorem}
Assume that $i+j\geq 2k$.  Then
\begin{align*}
& D^0+TD^0T=D+\phi_0, \\
& D^k+TD^{k}T=D^{k-1}\partial +\partial D^{k-1} +\phi_k
\end{align*}
as natural transformations 
$V_i\otimes W_j\to V_{i+j-k}\otimes W_{i+j-k}$.
\end{theorem}

\begin{proof}
The statement we need to prove is that
$A^k$ and $\phi_k$ induce the same natural transformation if
$i+j\geq 2k$. If $k=0$, this is (\ref{eq:Dwyercondition0}), which 
we have verified already.
Now we assume inductively that the statement is true for $k-1$, 
that is that $A^{k-1}$ and $\phi_{k-1}$ induce the same natural 
transformation if $i+j\geq 2k-2$. 

It follows that $A^{k-1}(d_0\otimes \id )$ and 
$\phi_{k-1}(d_0\otimes \id )$ 
induce the same transformation on $V_i\otimes W_j$ if $i+j\geq 2k-1$. 
But $\phi_{k-1}(d_0\otimes \id )$ is only non-trivial 
if $i=k$ and $j=k-1$, so  $A^{k-1}(d_0\otimes \id)$
induces the trivial natural transformation on $V_i\otimes W_j$ if
$i+j\geq 2k$. The same argument shows that
$A^{k-1}(\id\otimes d_0)$ is trivial on
$V_i\otimes W_j$ for $i+j\geq 2k$.

Because of of this and Lemma \ref{lemma:formula}, we get that
$A^k$ induces the same transformation as $\susp A^{k-1}$
on $V_i\otimes W_j$ for $i+j\geq 2k$. 
Now we use the induction assumption again, to see that
this transformation agrees with $\susp \phi_{k-1}=\phi_k$.
\end{proof}

\section{Appendix}

\begin{lemma}
\label{ztoz}
The formula in Definition \ref{deltadef} 
defines a map $\delta_i : Z_qR \to Z_{q+i}R$ when $2\leq i \leq q$.
\end{lemma}

\begin{proof}
Let $z\in Z_qR$. We must show that
$\delta_i (z)$ is a cycle ie.  
\[
d_j\delta_i (z) =
\sum_{ (\mu ,\nu ) \in \gshuf (q,i)}
d_j s_{\nu_i} \dots s_{\nu_1} (z) d_j s_{\mu_i} \dots s_{\mu_1} (z)=0,
\quad 0\leq j \leq q+i.
\] 
By one of the simplicial identities we have 
\[
d_is_j = 
\begin{cases}
s_{j-1}d_i &, i<j, \\
id &, i=j \text{ or } i=j+1, \\
s_jd_{i-1} &, i>j+1.
\end{cases}
\]

If $j\leq q-i$ we can commute the left $d_j$ all the way through 
the degeneracy maps to $z$ such that $d_j\delta_i (z)=0$.

If $j=q-i+1$ we find that all terms vanish except for those
with $\nu_1=q-i+1$, $\mu_1 = q-i$. For such a term we have
\[ 
d_js_{\nu_i} \dots s_{\nu_1} (z)d_js_{\mu_i} \dots s_{\mu_1} (z) =
s_{\nu_i -1} \dots s_{\nu_2 -1}(z)s_{\mu_i -1} \dots s_{\mu_2 -1}(z)
\]
so these cancels out in pairs.

If $j\geq q-i+2$ the non zero terms are those with $j \in \nu$,
$j-1\in \mu$ or $j-1\in \nu$, $j\in \mu$. 
For such a term we can interchange the corresponding $\mu_r$ 
with the corresponding $\nu_s$ and get a new element 
in $\gshuf (q,i)$ since $j-1 \geq q-i+1$. So these terms cancels
out in pairs.
\end{proof}

\begin{remark}
For $i=1$ we have that $\delta_1 (z)=s_q(z)s_{q-1}(z)$. Thus 
$d_j\delta_1 (z)=0$ for $j\neq q$ but $d_q\delta_1 (z)=z^2$.
\end{remark}

\end{document}